\documentclass[11pt]{amsart}
\usepackage{latexsym}
\usepackage{amsmath}
\usepackage{amssymb}

\newcommand{\eproof}{\mbox{\ }\hfill $\Box$ \par \vskip 10pt}

\newtheorem{Theorem}{Theorem}[section]
\newtheorem{lemma}[Theorem]{Lemma}
\newtheorem{prop}[Theorem]{Proposition}

\newtheorem{corol}[Theorem]{Corollary}

\def\cal{\mathcal}

\begin{document}

\title[Parabolic transmission eigenvalue-free regions]{Parabolic transmission eigenvalue-free regions in the degenerate isotropic case}

\author[G. Vodev]{Georgi Vodev}

\address {Universit\'e de Nantes, Laboratoire de Math\'ematiques Jean Leray, 2 rue de la Houssini\`ere, BP 92208, 44322 Nantes Cedex 03, France}
\email{Georgi.Vodev@univ-nantes.fr}

\date{}

\begin{abstract} We study the location of the transmission eigenvalues in the isotropic case when the restrictions of the refraction
indices on the boundary coincide. Under some natural conditions we show that there exist parabolic transmission eigenvalue-free regions.

\end{abstract} 

\maketitle

\setcounter{section}{0}
\section{Introduction and statement of results}

Let $\Omega\subset{\bf R}^d$, $d\ge 2$, be a bounded, connected domain with a $C^\infty$ smooth boundary $\Gamma=\partial\Omega$. 
A complex number $\lambda\neq 0$, ${\rm Re}\,\lambda\ge 0$, will be said to be a transmission eigenvalue if the following problem has a non-trivial solution:
$$\left\{
\begin{array}{lll}
\left(\Delta+\lambda^2 n_1(x)\right)u_1=0 &\mbox{in} &\Omega,\\
\left(\Delta+\lambda^2 n_2(x)\right)u_2=0 &\mbox{in} &\Omega,\\
u_1=u_2,\,\,\, \partial_\nu u_1=\partial_\nu u_2& \mbox{on}& \Gamma,
\end{array}
\right.
\eqno{(1.1)}
$$
where $\nu$ denotes the Euclidean unit inner normal to $\Gamma$, $n_j\in C^\infty(\overline\Omega)$, 
$j=1,2$ are strictly positive real-valued functions called refraction indices. In the non-degenerate isotropic case when
$$n_1(x)\neq n_2(x)\quad\mbox{on}\quad\Gamma\eqno{(1.2)}$$
it has been recently proved in \cite{kn:V3} that there are no transmission eigenvalues in the region 
$$\left\{\lambda\in{\bf C}:{\rm Re}\,\lambda\ge 0,\,\,|{\rm Im}\,\lambda|\ge C\right\}\eqno{(1.3)}$$
for some constant $C>0$. Moreover, it follows from the analysis
 in \cite{kn:LC} (see Section 4) that the eigenvalue-free region (1.3) is optimal and cannot be improved in general.
In the present paper we will consider the 
degenerate isotropic case when
$$n_1(x)\equiv n_2(x)\quad\mbox{on}\quad\Gamma.\eqno{(1.4)}$$
We suppose that there is an integer $j\ge 1$ such that
$$\partial_\nu^sn_1(x)\equiv \partial_\nu^sn_2(x)\quad\mbox{on}\quad\Gamma,\quad 0\le s\le j-1,\eqno{(1.5)}$$
$$\partial_\nu^jn_1(x)\neq \partial_\nu^jn_2(x)\quad\mbox{on}\quad\Gamma.\eqno{(1.6)}$$
It was proved in \cite{kn:CLM} (see Theorem 4.2) that in this case the eigenvalue-free region (1.3) is no longer valid. On the other hand, it follows from 
\cite{kn:LV} that under the conditions (1.5) and (1.6) there are no transmission eigenvalues in $|\arg\lambda|\ge\varepsilon$,
$|\lambda|\ge C_\varepsilon\gg 1$, 
$\forall\, 0<\varepsilon\ll 1$. Our goal in the present paper is to improve this result showing that in this case we have a much larger
parabolic eigenvalue-free region. Our main result is the following

\begin{Theorem} Under the conditions (1.5) and (1.6) there exists a constant $C>0$ such that there are no transmission eigenvalues in the region 
$$\left\{\lambda\in{\bf C}:{\rm Re}\,\lambda\ge 0,\,\,|{\rm Im}\,\lambda|\ge C\left({\rm Re}\,\lambda+1\right)^{1-\kappa_j
}\right\},\eqno{(1.7)}$$
where $\kappa_j=2(3j+2)^{-1}$.
\end{Theorem}
 
 To prove this theorem we make use of the semi-classical parametrix for the interior Dirichlet-to-Neumann (DN) map built in
 \cite{kn:V1}. It is proved in \cite{kn:V1} that for $|{\rm Im}\,\lambda|\ge \left({\rm Re}\,\lambda+1\right)^{1/2+\epsilon}$,
 $0<\epsilon\ll 1$, the DN map is an $h-\Psi$DO of class OP${\cal S}^1_{1/2-\epsilon}(\Gamma)$, where $0<h\ll 1$ is a semi-classical parameter
 such that $h\sim|\lambda|^{-1}$. A direct consequence of this fact is the existence of a transmission eigenvalue-free region of the form
 $$|{\rm Im}\,\lambda|\ge C_\epsilon\left({\rm Re}\,\lambda+1\right)^{1/2+\epsilon},\quad\forall\, 0<\epsilon\ll 1,\eqno{(1.8)}$$
under the condition (1.2). The most difficult part of the parametrix construction in \cite{kn:V1} is near the glancing region
(see Section 3 for the definition). Indeed, outside an arbitrary neighbourhood of the glancing region the parametrix construction in \cite{kn:V1}
works for $|{\rm Im}\,\lambda|\ge \left({\rm Re}\,\lambda+1\right)^{\epsilon}$ and the corresponding parametrix belongs to the class
OP${\cal S}^1_0(\Gamma)$. In other words, to improve the eigenvalue-free region (1.8) one has to improve the parametrix in the glancing region.
Such an improved parametrix has been built in \cite{kn:V2} for strictly concave domains and as a consequence (1.8) was improved to
$$|{\rm Im}\,\lambda|\ge C_\epsilon\left({\rm Re}\,\lambda+1\right)^{\epsilon},\quad\forall\, 0<\epsilon\ll 1,\eqno{(1.9)}$$
in this case. In fact, it turns out that to get larger eigenvalue-free regions under the condition (1.2) no parametrix construction
in the glancing region is needed. It suffices to show that the norm of the DN map microlocalized in a small neighbourhood of the glancing region
gets small if $|{\rm Im}\,\lambda|$ and ${\rm Re}\,\lambda$ are large. Indeed, this strategy has been implemented in \cite{kn:V3} to get the
optimal transmission eigenvalue-free region (1.3) for an arbitrary domain. In fact, the main point in the approach in \cite{kn:V3} is the
construction of a parametrix in the hyperbolic region valid for $1\ll C_\epsilon\le|{\rm Im}\,\lambda|\le \left({\rm Re}\,\lambda\right)^{1-\epsilon}$, ${\rm Re}\,\lambda\ge C'_\epsilon\gg 1$, $0<\epsilon\ll 1$.
The strategy of \cite{kn:V3}, however, does not work any more when we have the condition (1.4). In this case the parametrix
in the glancing region turns out to be essential to get eigenvalue-free regions like (1.7). In Section 3 we revisit the 
parametrix construction of \cite{kn:V1} and we study carefully the way in which it depends on the restriction on the boundary
of the normal derivatives of the refraction index (see Theorem 3.1). In Section 4 we improve Theorem 3.1. 
In Section 5 we show how Theorem 4.1 implies
Theorem 1.1. We also show that to improve (1.7) it suffices to improve the parametrix
in the glancing region, only (see Proposition 5.2).

 As in \cite{kn:PV1} one can study in this case the counting function
$N(r,C)=\#\{\lambda-{\rm trans.\, eig.}:\,C\le|\lambda|\le r\}$, where $r\gg C>0$. We have the following

\begin{corol}  Under the conditions of Theorem 1.1, there exists a constant $C>0$ such that the counting function of the transmission eigenvalues satisfies the asymptotics
$$N(r,C)=\tau r^d+{\cal O}_\varepsilon(r^{d-\kappa_j+\varepsilon}),\quad\forall\,0<\varepsilon\ll 1,\eqno{(1.10)}$$
where 
$$\tau=\frac{\omega_d}{(2\pi)^d}\int_\Omega \left(n_1(x)^{d/2}+n_2(x)^{d/2}\right)dx,$$
$\omega_d$ being the volume of the unit ball in ${\bf R}^d$.
\end{corol}

Note that the eigenvalue-free region (1.3) implies (1.10) with $\kappa_j$ replaced by $1$. 
Note also that asymptotics for the counting function $N(r,C)$ with remainder $o(r^d)$ have been previously 
obtained in \cite{kn:F}, \cite{kn:PS}, \cite{kn:R2} still under the condition (1.2).

\section{Basic properties of the $h-\Psi$DOs} 

In this section we will recall some basic properties of the $h-\Psi$DOs on a compact manifold without boundary.
Let $\Gamma$, ${\rm dim}\,\Gamma=d-1$, be as in the previous section and recall that given a 
symbol $a\in C^\infty(T^*\Gamma)$, the $h-\Psi$DO, ${\rm Op}_h(a)$, is defined as follows
$$\left({\rm Op}_h(a)f\right)(x')=(2\pi h)^{-d+1}\int_{T^*\Gamma}e^{-\frac{i}{h}\langle x'-y',\xi'\rangle}a(x',\xi')f(y')dy'd\xi'.$$
We have the following criteria of $L^2$- boundedness.

\begin{prop} Let the function $a$ satisfy the bounds
$$\left|\partial_{x'}^\alpha a(x',\xi')\right|\le C_\alpha,\quad\forall\,(x',\xi')\in T^*\Gamma, \eqno{(2.1)}$$
for all multi-indices $\alpha$. Then the operator ${\rm Op}_h(a)$ is bounded on $L^2(\Gamma)$ and satisfies
$$\left\|{\rm Op}_h(a)\right\|_{L^2(\Gamma)\to L^2(\Gamma)}\le C\sum_{|\alpha|\le d}C_\alpha\eqno{(2.2)}$$
with a constant $C>0$ independent of $h$ and $C_\alpha$.

Let the function $a$ satisfy the bounds
$$\left|\partial_{x'}^\alpha\partial_{\xi'}^\beta a(x',\xi')\right|\le C_{\alpha,\beta}h^{-(|\alpha|+|\beta|)/2},\quad\forall\,(x',\xi')\in T^*\Gamma, \eqno{(2.3)}$$
for all multi-indices $\alpha$ and $\beta$. Then the operator ${\rm Op}_h(a)$ is bounded on $L^2(\Gamma)$ and satisfies
$$\left\|{\rm Op}_h(a)\right\|_{L^2(\Gamma)\to L^2(\Gamma)}\le C\sum_{|\alpha|+|\beta|\le s_d}C_{\alpha,\beta}\eqno{(2.4)}$$
with a constant $C>0$ independent of $h$ and $C_{\alpha,\beta}$, where $s_d>0$ is an integer depending only on the dimension.
\end{prop}

Given $\ell\in{\bf R}$, $\delta_1,\delta_2\ge 0$ and a function $m>0$ on $T^*\Gamma$, we denote by
$S^\ell_{\delta_1,\delta_2}(m)$ the set of all functions $a\in C^\infty(T^*\Gamma)$ satisfying
$$\left|\partial_{x'}^\alpha\partial_{\xi'}^\beta a(x',\xi')\right|\le C_{\alpha,\beta}m^{\ell-\delta_1|\alpha|-
\delta_2|\beta|}$$
for all multi-indices $\alpha$ and $\beta$ with constants $C_{\alpha,\beta}>0$ independent of $m$. Given $k\in{\bf R}$, $0\le\delta<1/2$, we also denote by
${\cal S}_\delta^k$ the space of all symbols $a\in C^\infty(T^*\Gamma)$ satisfying
$$\left|\partial_{x'}^\alpha\partial_{\xi'}^\beta a(x',\xi')\right|\le C_{\alpha,\beta}h^{-\delta(|\alpha|+
|\beta|)}\langle\xi'\rangle^{k-|\beta|}$$
for all multi-indices $\alpha$ and $\beta$ with constants $C_{\alpha,\beta}>0$ independent of $h$. 
It is well-known that the $h-\Psi$DOs of class OP${\cal S}_\delta^k$ have nice calculus (e.g. see Section 7 of \cite{kn:DS}).
The next proposition is very usefull for inverting such operators depending on additional parameters (see also Proposition 2.2 of \cite{kn:V1}).

\begin{prop} Let $h^{\ell_\pm}a^\pm\in {\cal S}_\delta^{\pm k}$, $0\le\delta<1/2$, where $\ell_\pm\ge 0$ are some numbers.
Assume in addition that the functions $a^\pm$ satisfy
$$\left|\partial_{x'}^{\alpha_1}\partial_{\xi'}^{\beta_1}a^+(x',\xi')
\partial_{x'}^{\alpha_2}\partial_{\xi'}^{\beta_2}a^-(x',\xi')\right|\le\mu C_{\alpha_1,\beta_1,\alpha_2,\beta_2}h^{-(|\alpha_1|+|\beta_1|+|\alpha_2|+|\beta_2|)/2},\eqno{(2.5)}$$
$\forall\,(x',\xi')\in T^*\Gamma$,  
for all multi-indices $\alpha_1$, $\beta_1$, $\alpha_2$, $\beta_2$ such that $|\alpha_j|+|\beta_j|\ge 1$, $j=1,2$, with constants
$C_{\alpha_1,\beta_1,\alpha_2,\beta_2}>0$ independent of $h$ and $\mu$. Then we have
$$\left\|{\rm Op}_h(a^+){\rm Op}_h(a^-)-{\rm Op}_h(a^+a^-)\right\|_{L^2(\Gamma)\to L^2(\Gamma)}\le C(\mu+h)\eqno{(2.6)}$$
with a constant $C>0$ independent of $h$ and $\mu$.
\end{prop}

Given any real $s$, we define the semi-classical Sobolev norm by 
$$\|f\|_{H^s_h(\Gamma)}:=\|{\rm Op}_h(\langle\xi'\rangle^s)f\|_{L^2(\Gamma)}.$$
 Using the calculus of the $h-\Psi$DOs one can derive from
(2.4) the following

\begin{prop} Let $a\in {\cal S}_\delta^{-k}$, $0\le\delta<1/2$. Then, for every $s$, we have 
$${\rm Op}_h(a)={\cal O}_s(1): H_h^s(\Gamma)\to H_h^{s+k}(\Gamma).$$
\end{prop}

Proposition 2.2 implies the following 

 \begin{prop} Let $a^\pm\in {\cal S}_0^{\pm k}$. Then, for every $s$, we have  
 $${\rm Op}_h(a^+){\rm Op}_h(a^-)-{\rm Op}_h(a^+a^-)={\cal O}_s(h): H_h^s(\Gamma)\to H_h^{s}(\Gamma).$$
\end{prop}

\section{The parametrix construction revisited} 

In this section we will build a parametrix for the semi-classical Dirichlet-to-Neumann map following \cite{kn:V1}.
Note that in \cite{kn:V1} there is a gap due to a missing term in the transport equations (4.11), which however does not
affect the proof of the main results. Here we will correct this gap 
making some slight modifications.

Given $f\in H^{m+1}(\Gamma)$, let $u$ solve the equation 
$$\left\{
\begin{array}{lll}
\left(h^2\Delta+zn(x)\right)u=0 &\mbox{in} &\Omega,\\
u=f& \mbox{on}& \Gamma,
\end{array}
\right.
\eqno{(3.1)}
$$
where $n\in C^\infty(\overline\Omega)$ 
is a strictly positive function, $0<h\ll 1$ is a semi-classical parameter and $z\in Z_1\cup Z_2\cup Z_3$, where 
$Z_1=\{z\in {\bf C}:{\rm Re}\,z=1, 0<|{\rm Im}\,z|\le 1\}$, $Z_2=\{z\in {\bf C}:{\rm Re}\,z=-1, |{\rm Im}\,z|\le 1\}$, 
$Z_3=\{z\in {\bf C}:|{\rm Re}\,z|\le 1, |{\rm Im}\,z|=1\}$. Given $\varepsilon>0$ we also set 
$Z_1(\varepsilon)=\{z\in Z_1: h^{\varepsilon}\le|{\rm Im}\,z|\le 1\}$. We define 
the semi-classical Dirichlet-to-Neumann map 
$${\cal N}(h,z):H^{m+1}(\Gamma)\to H^{m}(\Gamma)$$
by
$${\cal N}(h,z)f:=-ih\partial_\nu u|_\Gamma$$
where $\nu$ denotes the Euclidean unit inner normal to $\Gamma$.
Given an integer $m\ge 0$, denote by $H_h^m(\Omega)$ the Sobolev space equipped with the semi-classical norm
$$\|v\|_{H_h^m(\Omega)}=\sum_{|\alpha|\le m}h^{|\alpha|}\left\|\partial_x^\alpha v\right\|_{L^2(\Omega)}.$$
We define similarly the Sobolev space $H_h^m(\Gamma)$. Note that this norm is equivalent to that one defined in Section 2.
Throughout this section we will use 
the normal coordinates $(x_1,x')$ with respect to the Euclidean metric near the boundary $\Gamma$, where
$0<x_1\ll 1$ denotes the Euclidean distance to the boundary and $x'$ are coordinates on $\Gamma$. We denote by
$\Delta_\Gamma$ the negative Laplace-Beltrami operator on $\Gamma$ equipped with the Riemannian metric induced by 
the Euclidean one in $\Omega$. Let $r_0(x',\xi')\ge 0$ be the principal symbol of $-\Delta_\Gamma$
written in the coordinates $(x',\xi')\in T^*\Gamma$. Since the function $n$ is smooth up to the boundary we can expand it as
$$n(x)=\sum_{k=0}^{N-1} x_1^kn_k(x')+ x_1^N{\cal M}_N(x)$$
for every integer $N\ge 1$, 
where $n_k=(k!)^{-1}\partial_\nu^k n|_\Gamma$, $n_0>0$, and ${\cal M}_N(x)$ is a real-valued smooth function. Set
$$\rho(x',\xi',z)=\sqrt{-r_0(x',\xi')+zn_0(x')},\quad {\rm Im}\,\rho>0.$$
The glancing region for the problem (3.1) is defined by
$$\Sigma:=\{(x',\xi')\in T^*\Gamma:r_\sharp(x',\xi')=1\},\quad r_\sharp=n_0^{-1}r_0.$$
Let $\phi\in C_0^\infty({\bf R})$, $0\le\phi\le 1$, $\phi(\sigma)=1$ for $|\sigma|\le 1$, 
$\phi(\sigma)=0$ for $|\sigma|\ge 2$, and set $\eta(x',\xi')=\phi(r_0(x',\xi')/\delta_0)$.
Clearly, taking $\delta_0>0$ small enough we can arrange that $|\rho|\ge C\langle r_0\rangle^{1/2}\ge C\langle\xi'\rangle$
on supp$(1-\eta)$. We also define the function $\chi(x',\xi')=\phi((1-r_\sharp(x',\xi'))/\delta_1)$, where $0<\delta_1\ll 1$ is independent of $h$
and $z$. Clearly, $\chi=1$ in a neighbourhood of $\Sigma$, $\chi=0$ outside another neighbourhood of $\Sigma$.

 We will say that a function $a\in C^\infty(T^*\Gamma)$ belongs to 
$S^{\ell_1}_{\delta_1,\delta_2}(m_1)+S^{\ell_2}_{\delta_3,\delta_4}(m_2)$ if $\eta a\in S^{\ell_1}_{\delta_1,\delta_2}(m_1)$
and $(1-\eta) a\in S^{\ell_2}_{\delta_3,\delta_4}(m_2)$. Given any integer $k$, it follows from Lemma 3.2 of \cite{kn:V1} that
$$\rho^k,\,|\rho|^k\in S^k_{2,2}(|\rho|)+S^k_{0,1}(|\rho|).\eqno{(3.2)}$$
In particular, (3.2) implies that
$$(1-\eta)\rho^k,\,(1-\eta)|\rho|^k\in{\cal S}_0^k.\eqno{(3.3)}$$
Since $\rho=i\sqrt{r_0}\left(1+{\cal O}(r_0^{-1})\right)$ as $r_0\to\infty$, it is easy to check that
$$(1-\eta)\rho^{k}-(1-\eta)(i\sqrt{r_0})^{k}\in {\cal S}_0^{k-2}\eqno{(3.4)}$$
for every integer $k$.
Since $|\rho|\ge C\sqrt{|{\rm Im}\,z|}$ for $z\in Z_1$, $(x',\xi')\in{\rm supp}\,\chi$ and $|\rho|\ge C>0$ for $z\in Z_2\cup Z_3$ 
or $z\in Z_1$, $(x',\xi')\in{\rm supp}\,(1-\chi)$ (see Lemma 3.1 of \cite{kn:V1}), it also follows from (3.2) that 
$$(1-\chi)\rho^k,\,(1-\chi)|\rho|^k\in{\cal S}_0^k,\eqno{(3.5)}$$
$$h^{\frac{k_-}{4}}\chi\rho^k,\,h^{\frac{k_-}{4}}\chi|\rho|^k\in{\cal S}_{1/2-\epsilon}^{-N},\quad z\in Z_{1}(1/2-\epsilon),\eqno{(3.6)}$$
$$\chi\rho^k,\,\chi|\rho|^k\in{\cal S}_0^{-N},\quad z\in Z_2\cup Z_3,\eqno{(3.7)}$$
for every integer $N\ge 0$ and $0<\epsilon\ll 1$, where $k_-=0$ if $k\ge 0$,
$k_-=|k|$ if $k<0$. Our goal in this section is to prove the following

\begin{Theorem} Let $z\in Z_1(1/2-\epsilon)$, $0<\epsilon\ll 1$. Then, for every integer $s\ge 0$ there is a function $b_s\in {\cal S}_{1/2-\epsilon}^0$ independent of all $n_k$ with $k\ge s$ such that
$$\left\|{\cal N}(h,z)-{\rm Op}_h\left(\rho+b_s+c_sh^s\rho^{-s-1}zn_s\right)\right\|_{L^2(\Gamma)\to H_h^{s+1}(\Gamma)}\le 
C_sh^{s+1}|{\rm Im}\,z|^{-2s-3/2}\eqno{(3.8)}$$
where $c_s=0$ if $s=0$, and $c_s=-is!(-2i)^{-s-1}$ for $s\ge 1$.
If $z\in Z_2\cup Z_3$, then (3.8) holds with $|{\rm Im}\,z|$ replaced by $1$.
 Moreover, for $z\in Z_1(1-\epsilon)$ we have
$$\left\|{\cal N}(h,z){\rm Op}_h(1-\chi)-{\rm Op}_h\left(\rho(1-\chi)+\widetilde b_s+c_sh^s(1-\chi)\rho^{-s-1}zn_s\right)\right\|_{L^2(\Gamma)\to H_h^{s+1}(\Gamma)}\le 
C_sh^{s+1}\eqno{(3.9)}$$
where the function $\widetilde b_s\in {\cal S}_0^0$ is independent of all $n_k$ with $k\ge s$.
\end{Theorem}

{\it Proof.} We will recall the parametrix construction in \cite{kn:V1}. 
We will proceed locally and then we will use partition of the unity to get the global parametrix. Fix a point 
$x^0\in\Gamma$ and let ${\cal U}_0\subset\Gamma$ be a small open neighbourhood of $x^0$. Let $(x_1,x')$, $x_1>0$, $x'\in 
{\cal U}_0$, be the normal coordinates. In these coordinates the Laplacian can be written as follows
$$\Delta=\partial_{x_1}^2+r(x,\partial_{x'})+q(x,\partial_x)$$
where $r(x,\xi')=\langle R(x)\xi',\xi'\rangle\ge 0$, $R=(R_{ij})$ being a symmetric $(d-1)\times(d-1)$ matrix-valued function with smooth
real-valued entries, $q(x,\xi)=\langle q(x),\xi\rangle=q^\sharp(x)\xi_1+\langle q^\flat(x),\xi'\rangle$, $q^\sharp$ and $q^\flat$ being smooth functions. We can expand them as follows
$$R(x)=\sum_{k=0}^{N-1} x_1^kR_k(x')+x_1^N{\cal R}_N(x),$$
$$q^\sharp(x)=\sum_{k=0}^{N-1} x_1^kq^\sharp_k(x')+x_1^N{\cal Q}^\sharp_N(x),$$ 
$$q^\flat(x)=\sum_{k=0}^{N-1} x_1^kq^\flat_k(x')+x_1^N{\cal Q}^\flat_N(x),$$
for every integer $N\ge 1$. Clearly, $r_0(x',\xi')=r(0,x',\xi')=\langle R_0(x')\xi',\xi'\rangle$.

Take a function $\psi^0\in C_0^\infty({\cal U}_0)$. In what follows $\psi$ will denote either the function $\psi^0$ or
the function $\psi^0(1-\chi)$. 
Following \cite{kn:V1}, we will construct a parametrix $\widetilde u_\psi$ of the solution of
(3.1) with $\widetilde u_\psi|_{x_1=0}={\rm Op}_h(\psi)f$ in the form  
$$\widetilde u_\psi=(2\pi h)^{-d+1}\int\int e^{\frac{i}{h}(\langle y',\xi'\rangle+\varphi(x,\xi',z))}\Phi_\delta(x,\xi',z)a(x,\xi',h,z)f(y')
d\xi'dy'$$
where $\Phi_\delta=\phi(x_1/\delta)\phi(x_1/\rho_1\delta)$, with $\rho_1=|\rho|^3$ if $z\in Z_{1}(1/2-\epsilon)$, $\psi=\psi^0$,
and $\rho_1=1$ if 
$z\in Z_2\cup Z_3$ or $z\in Z_{1}(1-\epsilon)$, $\psi=\psi^0(1-\chi)$. Here 
$0<\delta\ll 1$ is a parameter independent of $h$ and $z$ to be fixed later on. 
The phase $\varphi$ is complex-valued such that $\varphi|_{x_1=0}=-\langle x',\xi'\rangle$ and satisfies the eikonal equation
mod ${\cal O}(x_1^N)$:
$$\left(\partial_{x_1}\varphi\right)^2+\left\langle R(x)\nabla_{x'}\varphi,\nabla_{x'}\varphi\right\rangle-zn(x)=
x_1^N\Psi_N\eqno{(3.10)}$$
where $N\gg 1$ is an arbitrary integer and the function $\Psi_N$ is smooth up to the boundary $x_1=0$. 
It is shown in \cite{kn:V1}, Section 4, that the equation (3.10) has a smooth solution
of the form
$$\varphi=\sum_{k=0}^{N}x_1^k\varphi_k(x',\xi',z),\quad \varphi_0=-\langle x',\xi'\rangle,$$
safisfying 
$$\partial_{x_1}\varphi|_{x_1=0}=\varphi_1=\rho.\eqno{(3.11)}$$ 
More generally, the functions $\varphi_k$ satisfy the relations
$$\sum_{k+j=K}(k+1)(j+1)\varphi_{k+1}\varphi_{j+1}+\sum_{k+j+\ell=K}\left\langle R_\ell\nabla_{x'}\varphi_k,\nabla_{x'}
\varphi_j\right\rangle-zn_K=0\eqno{(3.12)}$$ 
for every integer $0\le K\le N-1$. Then equation (3.10) is satisfied with
$$\Psi_N=\left\langle{\cal R}_N(x)\nabla_{x'}\varphi,\nabla_{x'}\varphi\right\rangle-z{\cal M}_N(x)$$
$$+\sum_{k+j\ge N}x_1^{k+j-N}(k+1)(j+1)\varphi_{k+1}\varphi_{j+1}+\sum_{k+j+\ell\ge N}x_1^{k+j+\ell-N}\left\langle R_\ell\nabla_{x'}\varphi_k,\nabla_{x'}\varphi_j\right\rangle$$
where $\varphi_\nu=0$ for $\nu\ge N+1$ so that the above sums are finite. 
Using (3.12) one can prove by induction the following lemma (see Lemma 4.1 of \cite{kn:V1}).

\begin{lemma} We have
$$\varphi_k\in S^{4-3k}_{2,2}(|\rho|)+S^1_{0,1}(|\rho|),\quad 1\le k\le N,\eqno{(3.13)}$$
$$\partial_{x_1}^k\Psi_N\in S^{2-3N-3k}_{2,2}(|\rho|)+S^2_{0,1}(|\rho|),\quad k\ge 0,\eqno{(3.14)}$$
uniformly in $z$ and $0\le x_1\le 2\delta\min\{1,|\rho|^3\}$. Moreover, if $\delta>0$ is small enough, independent of
$\rho$, we have
$${\rm Im}\,\varphi\ge x_1{\rm Im}\,\rho/2\quad\mbox{for}\quad 0\le x_1\le 2\delta\min\{1,|\rho|^3\}.\eqno{(3.15)}$$
\end{lemma}

One can also easily prove by induction the following

\begin{lemma} For every integer $k\ge 1$ the functions $\varphi_k$ and $\varphi_{k+1}-\frac{zn_k}{2(k+1)\rho}$ are independent of all $n_\ell$ with $\ell\ge k$. 
\end{lemma}

It follows from (3.13) that $(1-\eta)\varphi_k\in {\cal S}_0^1$ for all $k$. Define now the functions $\widetilde \varphi_k$
independent of all $n_\ell$, $\ell\ge 0$, satisfying the relations
$$\sum_{k+j=K}(k+1)(j+1)\widetilde\varphi_{k+1}\widetilde\varphi_{j+1}+\sum_{k+j+\ell=K}\left\langle R_\ell\nabla_{x'}
\widetilde\varphi_k,\nabla_{x'}
\widetilde\varphi_j\right\rangle=0,\eqno{(3.16)}$$ 
$1\le k\le N-1$, and 
$\widetilde\varphi_0=-\langle x',\xi'\rangle$, $\widetilde\varphi_1=i\sqrt{r_0}$. Using (3.4) together with (3.12) and (3.16), one can easily prove by induction the following

\begin{lemma} For every integer $k\ge 1$, we have $(1-\eta)(\varphi_k-\widetilde\varphi_k)\in {\cal S}_0^{-1}$.
\end{lemma}

The amplitude $a$ is of the form
$$a=\sum_{j=0}^{N-1} h^{j}a_j(x,\xi',z)$$
where the functions $a_j$ satisfy the transport equations mod ${\cal O}(x_1^N)$:
$$2i\partial_{x_1}\varphi\partial_{x_1}a_j+
2i\left\langle R(x)\nabla_{x'}\varphi,\nabla_{x'}a_j\right\rangle+i\left(\Delta\varphi\right)a_j+\Delta a_{j-1}=x_1^N A_N^{(j)},\quad 0\le j\le N-1,\eqno{(3.17)}$$
$a_0|_{x_1=0}=\psi$, $a_j|_{x_1=0}=0$ for $j\ge 1$, where $a_{-1}=0$ and the functions $A_N^{(j)}$ are
smooth up to the boundary $x_1=0$. We will be looking for the solutions to (3.17) in the form
$$a_j=\sum_{k=0}^{N}x_1^ka_{k,j}(x',\xi',z).$$
We can write
$$\Delta\varphi=\sum_{k=0}^{N-1}x_1^k\varphi_k^\Delta+x_1^NE_N(x)$$
with
$$\varphi_k^\Delta=(k+1)(k+2)\varphi_{k+2}+\sum_{\ell+\nu=k}\left(\langle R_\ell\nabla_{x'}\cdot,\nabla_{x'}\varphi_\nu\rangle+
q_\ell^\sharp(\nu+1)\varphi_{\nu+1}+\langle q_\ell^\flat,\nabla_{x'}\varphi_\nu\rangle\right),$$
$$E_N=\langle{\cal R}_N\nabla_{x'}\cdot,\nabla_{x'}\varphi\rangle+{\cal Q}_N^\sharp\partial_{x_1}\varphi+\langle {\cal Q}_N^\flat,
\nabla_{x'}\varphi\rangle$$
 $$+\sum_{\ell+\nu\ge N}x_1^{\ell+\nu-N}\left(\langle R_\ell\nabla_{x'}\cdot,\nabla_{x'}\varphi_\nu\rangle+
q_\ell^\sharp(\nu+1)\varphi_{\nu+1}+\langle q_\ell^\flat,\nabla_{x'}\varphi_\nu\rangle\right),$$
where $\varphi_\nu=0$ for $\nu\ge N+1$. 
Similarly 
$$\Delta a_{j-1}=\sum_{k=0}^{N-1}x_1^ka_{k,j-1}^\Delta+x_1^NF^{(j-1)}_N(x)$$
with
$$a_{k,j-1}^\Delta=(k+1)(k+2)a_{k+2,j-1}$$ $$+\sum_{\ell+\nu=k}\left(\langle R_\ell\nabla_{x'}\cdot,\nabla_{x'}a_{\nu,j-1}\rangle+
q_\ell^\sharp(\nu+1)a_{\nu+1,j-1}+\langle q_\ell^\flat,\nabla_{x'}a_{\nu,j-1}\rangle\right),$$
$$F^{(j-1)}_N=\langle{\cal R}_N\nabla_{x'}\cdot,\nabla_{x'}a_{j-1}\rangle+{\cal Q}_N^\sharp\partial_{x_1}a_{j-1}+\langle {\cal Q}_N^\flat,
\nabla_{x'}a_{j-1}\rangle$$
$$+\sum_{\ell+\nu\ge N}x_1^{\ell+\nu-N}\left(\langle R_\ell\nabla_{x'}\cdot,\nabla_{x'}a_{\nu,j-1}\rangle+
q_\ell^\sharp(\nu+1)a_{\nu+1,j-1}+\langle q_\ell^\flat,\nabla_{x'}a_{\nu,j-1}\rangle\right),$$
where $a_{\nu,j-1}=0$ for $\nu\ge N+1$. 
We also have
$$\left(\Delta\varphi\right)a_j=\sum_{k=0}^{N-1}x_1^k\sum_{k_1+k_2=k}\varphi_{k_1}^\Delta a_{k_2,j}+x_1^N{\cal E}_N^{(j)}$$
with
$${\cal E}_N^{(j)}=E_Na_j+\sum_{k_1+k_2\ge N}x_1^{k_1+k_2-N}\varphi_{k_1}^\Delta a_{k_2,j},$$
$$\partial_{x_1}\varphi\partial_{x_1}a_j=\sum_{k=0}^{N-1}x_1^k\sum_{k_1+k_2=k}(k_1+1)(k_2+1)\varphi_{k_1+1} a_{k_2+1,j}+x_1^N{\cal F}_N^{(j)}$$
with
$${\cal F}_N^{(j)}=\sum_{k_1+k_2\ge N}x_1^{k_1+k_2-N}(k_1+1)(k_2+1)\varphi_{k_1+1} a_{k_2+1,j},$$
$$\left\langle R(x)\nabla_{x'}\varphi,\nabla_{x'}a_j\right\rangle=\sum_{k=0}^{N-1}x_1^k\sum_{k_1+k_2+k_3=k}
\left\langle R_{k_1}\nabla_{x'}\varphi_{k_2},\nabla_{x'}a_{k_3,j}\right\rangle+x_1^N{\cal G}_N^{(j)}$$
with 
$${\cal G}_N^{(j)}=\left\langle{\cal R}_N(x)\nabla_{x'}\varphi,\nabla_{x'}a_j\right\rangle+\sum_{k_1+k_2+k_3\ge N}x_1^{k_1+k_2+k_3-N}
\left\langle R_{k_1}\nabla_{x'}\varphi_{k_2},\nabla_{x'}a_{k_3,j}\right\rangle,$$
where $\varphi_\nu=0$, $a_{\nu,j}=0$ for $\nu\ge N+1$ so that the above sums are finite. 
Inserting the above identities into equation (3.17) and comparing the coefficients of all powers $x_1^k$, $0\le k\le N-1$, we get that the functions $a_{k,j}$ must satisfy the relations
$$\sum_{k_1+k_2=k}2i(k_1+1)(k_2+1)\varphi_{k_1+1} a_{k_2+1,j}+\sum_{k_1+k_2+k_3=k}
2i\left\langle R_{k_1}\nabla_{x'}\varphi_{k_2},\nabla_{x'}a_{k_3,j}\right\rangle$$ $$
+\sum_{k_1+k_2=k}i\varphi_{k_1}^\Delta a_{k_2,j}=-a_{k,j-1}^\Delta,\quad\mbox{for}\quad  0\le k\le N-1,\quad 0\le j\le N-1,\eqno{(3.18)}$$
and $a_{0,0}=\psi$, $a_{0,j}=0$, $j\ge 1$, $a_{k,-1}=0$, $k\ge 0$. Then equation (3.17) is satisfied with
$$A_N^{(j)}=
2i{\cal F}_N^{(j)}+2i{\cal G}_N^{(j)}+i{\cal E}_N^{(j)}+F_N^{(j-1)}.$$
Let us calculate $a_{1,0}$.
By (3.18) with $j=0$, $k=0$, we get
$$a_{1,0}=-\varphi_1^{-1}\langle B_0\xi',\nabla_{x'}\psi\rangle-(\varphi_1^{-1}\varphi_2+2^{-1}q^\sharp_0-(2\varphi_1)^{-1}\langle q^\flat_0(x'),\xi'\rangle)\psi.\eqno{(3.19)}$$
On the other hand, by (3.12) with $K=1$ we get
$$\varphi_2=-(2\rho)^{-1}\langle B_0\xi',\nabla_{x'}\rho\rangle-(4\rho)^{-1}\langle B_1\xi',\xi'\rangle +z(4\rho)^{-1}n_1.$$
Using the identity 
$$2\rho\nabla_{x'}\rho=-\nabla_{x'}r_0+z\nabla_{x'}n_0$$
we can write $\varphi_2$ in the form
$$\varphi_2=(2\rho)^{-2}\langle B_0\xi',\nabla_{x'}r_0\rangle-(4\rho)^{-1}\langle B_1\xi',\xi'\rangle$$ $$-z(2\rho)^{-2}\langle B_0\xi',\nabla_{x'}n_0\rangle+z(4\rho)^{-1}n_1.\eqno{(3.20)}$$
By (3.19) and (3.20),
$$a_{1,0}=-\rho^{-1}\langle B_0\xi',\nabla_{x'}\psi\rangle-2^{-1}\psi q_0^\sharp+(2\rho)^{-1}\psi\langle q_0^\flat(x'),\xi'\rangle$$
 $$-4^{-1}\rho^{-3}\psi\langle B_0\xi',\nabla_{x'}r_0\rangle+4^{-1}\rho^{-2}\psi\langle B_1\xi',\xi'\rangle$$ $$+z4^{-1}\rho^{-3}\psi\langle B_0\xi',\nabla_{x'}n_0\rangle-z4^{-1}\rho^{-2}\psi n_1.\eqno{(3.21)}$$
By (3.2) and (3.21) we conclude
 $$a_{1,0}\in S^{-3}_{2,2}(|\rho|)+S^0_{0,1}(|\rho|).\eqno{(3.22)}$$
 The next lemma follows from Lemma 3.2 and (3.22) together with equations (3.18) and can be proved in 
 the same way as Lemma 4.2 of \cite{kn:V1}. We will sketch the proof.
 
 \begin{lemma} We have
 $$a_{k,j}\in S^{-3k-4j}_{2,2}(|\rho|)+S^{-j}_{0,1}(|\rho|),\quad\mbox{for}\quad k\ge 1,\,j\ge 0,\eqno{(3.23)}$$
 $$\partial_{x_1}^kA_N^{(j)}\in S^{-3N-3k-4j-2}_{2,2}(|\rho|)+S^{1-j}_{0,1}(|\rho|),\quad\mbox{for}\quad k\ge 0,\,j\ge 0,\eqno{(3.24)}$$
 uniformly in $z$ and $0\le x_1\le 2\delta\min\{1,|\rho|^3\}$.
 \end{lemma}
 
 {\it Proof.} Recall that $\nabla_{x'}\varphi_0=-\xi'$. By (3.13) we have
 $$\nabla_{x'}\varphi_k\in S^{2-3k}_{2,2}(|\rho|)+S^1_{0,1}(|\rho|),\quad k\ge 1,$$
 $$\varphi_k^\Delta\in S^{-2-3k}_{2,2}(|\rho|)+S^1_{0,1}(|\rho|),\quad k\ge 0.$$
 We will prove (3.23) by induction. In view of (3.22) we have (3.23) with $k=1$, $j=0$.
Suppose now that (3.23) is true for all $j\le J-1$ and all $k\ge 1$, and for $j=J$ and $k\le K$. We have to show that it is true
for $j=J$ and $k=K+1$. To this end, we will use equation (3.18) with $j=J$ and $k=K$. Indeed, 
the LHS is equal to $2i(K+1)\rho a_{K+1,J}$ modulo $S^{-3K-4J-2}_{2,2}(|\rho|)+S^{-J+1}_{0,1}(|\rho|)$, while the RHS belongs to
$S^{-3K-4J-2}_{2,2}(|\rho|)+S^{-J+1}_{0,1}(|\rho|)$. In other words, $\rho a_{K+1,J}$ belongs to
$S^{-3K-4J-2}_{2,2}(|\rho|)+S^{-J+1}_{0,1}(|\rho|)$. This implies that $a_{K+1,J}$ belongs to
$S^{-3K-4J-3}_{2,2}(|\rho|)+S^{-J}_{0,1}(|\rho|)$, as desired. Furthermore, (3.24) follows from (3.13) and (3.23) since 
the functions $A_N^{(j)}$ are expressed in terms of $\varphi_k$ and $a_{k,j}$.
One needs the simple observation that
$$a\in S_{2,2}^{\ell_1}(|\rho|)+S_{0,1}^{\ell_2}(|\rho|)$$
implies
$$x_1^k a\in S_{2,2}^{\ell_1+3k}(|\rho|)+S_{0,1}^{\ell_2}(|\rho|).$$
 \eproof
 
Using Lemma 3.3 we will prove the following

\begin{lemma} For all $k\ge 1,\,j\ge 0$, the function
$$a_{k,j}-\frac{(k+j)!}{k!}\frac{z\psi n_{k+j}}{(-2i\rho)^{j+2}}$$
is independent of all $n_\ell$ with $\ell\ge k+j$.
\end{lemma}

{\it Proof.} It follows from Lemma 3.3 that the function
$$\varphi_k^\Delta-(2\rho)^{-1}(k+1)zn_{k+1}$$
is independent of all $n_\ell$ with $\ell\ge k+1$. We will first prove the assertion for $j=0$ and all $k\ge 1$
by induction in $k$. In view of (3.21) it is true for $k=1$. Suppose it is true for all integers $k\le K$ with some integer $K\ge 1$.
We will prove it for $k=K+1$. To this end, we will use equation (3.18) with $j=0$ and $k=K$. Since the RHS is zero, we get that the function
$$2i(K+1)\rho a_{K+1,0}+i\varphi_K^\Delta\psi$$
is independent of all $n_\ell$ with $\ell\ge K+1$. Hence, so is the function
$$a_{K+1,0}+(2\rho)^{-2}(K+1)z\psi n_{K+1}$$
as desired. We will now prove the assertion for all $k\ge 1$, $j\ge 0$ by induction in $j$. Suppose it is true for $j\le J$ and all $k\ge 1$
with some integer $J\ge 1$. We will prove it for $j=J+1$ and all $k\ge 1$. To this end, we will use equation (3.18) with $j=J+1$
and $k$ replaced by $k-1$, $k\ge 1$. We have that, modulo functions independent of all $n_\ell$ with $\ell\ge k+J+1$, the LHS is equal
to $2ik\rho a_{k,J+1}$, while the RHS is equal to $-k(k+1)a_{k+1,J}$. Hence the function
$$a_{k,J+1}+(2i\rho)^{-1}(k+1)a_{k+1,J}$$
is independent of all $n_\ell$ with $\ell\ge k+J+1$, which clearly implies the desired assertion.
\eproof 

It follows from (3.23) that $(1-\eta)a_{k,j}\in {\cal S}_0^{-j}$ for all $k\ge 1$, $j\ge 0$. 
Define now the functions $\widetilde a_{k,j}$
independent of all $n_\ell$, $\ell\ge 0$, satisfying the relations
$$\sum_{k_1+k_2=k}2i(k_1+1)(k_2+1)\widetilde\varphi_{k_1+1}\widetilde a_{k_2+1,j}+\sum_{k_1+k_2+k_3=k}
2i\left\langle R_{k_1}\nabla_{x'}\widetilde\varphi_{k_2},\nabla_{x'}\widetilde a_{k_3,j}\right\rangle$$ $$
+\sum_{k_1+k_2=k}i\widetilde\varphi_{k_1}^\Delta\widetilde a_{k_2,j}=-\widetilde a_{k,j-1}^\Delta,\eqno{(3.25)}$$
and $\widetilde a_{0,0}=\psi$, $\widetilde a_{0,j}=0$, $j\ge 1$, $\widetilde a_{k,-1}=0$, $k\ge 0$, where
$\widetilde\varphi_{k}^\Delta$ is defined by replacing in the definition of $\varphi_{k}^\Delta$ all functions $\varphi_j$
by $\widetilde\varphi_j$.
Using Lemma 3.4 we will prove the following

\begin{lemma} For all $k\ge 1,\,j\ge 0$, we have $(1-\eta)(a_{k,j}-\widetilde a_{k,j})\in {\cal S}_0^{-j-1}$.
\end{lemma}

{\it Proof.}  By Lemma 3.4 together with (3.18) and (3.25) we obtain that the relations
$$\sum_{k_1+k_2=k}2i(k_1+1)(k_2+1)(1-\eta)\widetilde\varphi_{k_1+1}(a_{k_2+1,j}-\widetilde a_{k_2+1,j})$$ $$+\sum_{k_1+k_2+k_3=k}
2i(1-\eta)\left\langle R_{k_1}\nabla_{x'}\widetilde\varphi_{k_2},\nabla_{x'}(a_{k_3,j}-\widetilde a_{k_3,j})\right\rangle$$ $$
+\sum_{k_1+k_2=k}i(1-\eta)\widetilde\varphi_{k_1}^\Delta(a_{k_2,j}-\widetilde a_{k_2,j})=-(1-\eta)(a_{k,j-1}^\Delta-\widetilde a_{k,j-1}^\Delta)\eqno{(3.26)}$$
are satisfied modulo ${\cal S}_0^{-j-1}$. We will proceed by induction. 
Suppose now that the assertion is true for all $j\le J-1$ and all $k\ge 1$, and for $j=J$ and $k\le K$. This implies that the LHS of 
(3.26) with $k=K$ and $j=J$ is equal to $2i(K+1)(1-\eta)\widetilde\varphi_1(a_{K+1,J}-\widetilde a_{K+1,J})$ modulo 
${\cal S}_0^{-J}$, while the RHS belongs to ${\cal S}_0^{-J}$. Hence, $(1-\eta)(a_{K+1,J}-\widetilde a_{K+1,J})$
belongs to ${\cal S}_0^{-J-1}$, as desired.
\eproof
 
In view of (3.11) we have $$-ih\partial_{x_1}\widetilde u_\psi|_{x_1=0}={\cal T}_\psi(h,z)f={\rm Op}_h(\tau_\psi)f$$ where
$$\tau_\psi=a\frac{\partial\varphi}{\partial x_1}|_{x_1=0}-ih\frac{\partial a}{\partial x_1}|_{x_1=0}=\rho\psi-i\sum_{j=0}^{N-1}h^{j+1}a_{1,j}.$$

\begin{lemma} For every integer $m\ge 0$ there are $N_m>1$ and $\ell_m>0$ such that for all
$N\ge N_m$ we have the estimate
$$\left\|{\cal N}(h,z){\rm Op}_h(\psi)-{\cal T}_\psi(h,z)\right\|_{L^2(\Gamma)\to H_h^m(\Gamma)}\le C_{N,m}h^{\epsilon N-\ell_m}\eqno{(3.27)}$$
if $\psi=\psi^0$, $z\in Z_1(1/2-\epsilon)$, or $\psi=\psi^0(1-\chi)$, $z\in Z_1(1-\epsilon)$. 
If $\psi=\psi^0$, $z\in Z_2\cup Z_3$, then (3.27) holds with $\epsilon$ replaced by $1$.
\end{lemma}

{\it Proof.} Denote by $G_D$ the Dirichlet self-adjoint realization of the operator $-n^{-1}\Delta$ on the Hilbert space 
$L^2(\Omega;n(x)dx)$. It is easy to see that
$$(h^2G_D-z)^{-1}={\cal O}\left(\theta(z)^{-1}\right):L^2(\Omega)\to L^2(\Omega)$$
where $\theta(z)=|{\rm Im}\,z|$ if $z\in Z_1$, $\theta(z)=1$ if $z\in Z_2\cup Z_3$.
Clearly, under the conditions of Lemma 3.8, we have $h<\theta(z)\le 1$. The above bound together with the coercivity of $G_D$ imply
$$(h^2G_D-z)^{-1}={\cal O}_s\left(\theta(z)^{-1}\right):H^s(\Omega)\to H^s(\Omega)\eqno{(3.28)}$$
for every integer $s\ge 0$. 
We also have the identity
$${\cal N}(h,z){\rm Op}_h(\psi)f-{\cal T}_\psi(h,z)f=-ih\gamma\partial_\nu\left((h^2G_D-z)^{-1}V\right)\eqno{(3.29)}$$
where $\gamma$ denotes the restriction on $\Gamma$, and 
$$V=(h^2\Delta+zn)\widetilde u_\psi$$
$$={\cal K}(h,z)f=(2\pi h)^{-d+1}\int\int e^{\frac{i}{h}(\langle y',\xi'\rangle+\varphi(x,\xi',z))}K(x,\xi',h,z)f(y')
d\xi'dy',$$
where $K=K_1+K_2$ with
$$K_1=\left[h^2\Delta,\Phi_\delta\right]a,\quad K_2=\left(x_1^NA_N+h^NB_N\right)\Phi_\delta,$$
$$A_N=\Psi_Na+\sum_{j=0}^{N-1}h^{j+1}A_N^{(j)},\quad B_N=\Delta a_{N-1}=\sum_{k=0}^{N-1} x_1^ka_{k,N-1}^\Delta+x_1^NF_N^{(N-1)}.$$
By the trace theorem we get from (3.28) and (3.29),
$$\left\|{\cal N}(h,z){\rm Op}_h(\psi)f-{\cal T}_\psi(h,z)f\right\|_{H_h^m(\Gamma)}\le {\cal O}(h^{-1})\left\|(h^2G_D-z)^{-1}V\right\|_{H_h^{m+1}(\Omega)}$$
$$\le {\cal O}\left((h\theta(z))^{-1}\right)\left\|V\right\|_{H_h^{m+1}(\Omega)}
\le {\cal O}\left(h^{-2}\right)\left\|V\right\|_{H_h^{m+1}(\Omega)}.\eqno{(3.30)}$$
To bound the norm of $V$ we need to bound the kernel of the operator 
$${\cal K}_\alpha:=\partial_x^\alpha{\cal K}(h,z):L^2(\Gamma)\to L^2(\Omega).$$ 
By Lemma 3.1 of \cite{kn:V1} we have
$${\rm Im}\,\rho\ge \frac{|{\rm Im}\,z|}{2|\rho|}\quad \mbox{on}\quad{\rm supp}\,\eta,\,z\in Z_1,$$
$${\rm Im}\,\rho\ge C\langle\xi'\rangle \quad \mbox{for}\quad z\in Z_2\cup Z_3\quad\mbox{and on}\quad{\rm supp}\,(1-\eta),\,z\in Z_1,$$
where $C>0$ is some constant. Hence, by (3.15), for $0\le x_1\le 2\delta\min\{1,|\rho|^3\}$ we have
$$x_1^N\left|e^{i\varphi/h}\right|\le x_1^N e^{-{\rm Im}\,\varphi/h}\le x_1^N e^{-x_1{\rm Im}\,\rho/2h}
\le C_N\left(\frac{h}{{\rm Im}\,\rho}\right)^N$$
$$\le\left\{
\begin{array}{ll}
C_N\left(\frac{h|\rho|}{|{\rm Im}\,z|}\right)^N, &z\in Z_1,\,(x',\xi')\in{\rm supp}\,\eta,\\
C_N\left(\frac{h}{\langle\xi'\rangle }\right)^N, &\mbox{otherwise}.
\end{array}
\right.
\eqno{(3.31)}
$$
On the other hand, by Lemmas 3.2 and 3.5, for $|\rho|^4\ge h$ and $0\le x_1\le 2\delta\min\{1,|\rho|^3\}$ we have
$$\left|\partial_x^\alpha A_N\right|
\le\left\{
\begin{array}{ll}
C_{\alpha,N}h^{-\ell_\alpha}|\rho|^{-3N}, &\mbox{on}\quad{\rm supp}\,\eta,\\
C_{\alpha,N}\langle\xi'\rangle^2, &\mbox{on}\quad{\rm supp}\,(1-\eta),
\end{array}
\right.
\eqno{(3.32)}
$$
$$\left|\partial_x^\alpha B_N\right|
\le\left\{
\begin{array}{ll}
C_{\alpha,N}h^{-\ell_\alpha}|\rho|^{-4N}, &\mbox{on}\quad{\rm supp}\,\eta,\\
C_{\alpha,N}\langle\xi'\rangle^{-N+1}, &\mbox{on}\quad{\rm supp}\,(1-\eta),
\end{array}
\right.
\eqno{(3.33)}
$$
for every multi-index $\alpha$ with some $\ell_\alpha>0$ independent of $N$. By (3.31), (3.32) and (3.33), using that
$|\rho|^2\ge C|{\rm Im}\,z|$, $C>0$, on supp$\,\eta$, we conclude
$$\left|\partial_x^\alpha \left(e^{i\varphi/h}K_2\right)\right|\le
C_{\alpha,N}h^{-\ell_\alpha}\left(\frac{h}{|\rho|^2|{\rm Im}\,z|}+\frac{h}{|\rho|^4}\right)^N\le
C_{\alpha,N}h^{-\ell_\alpha}\left(\frac{h}{|{\rm Im}\,z|^2}\right)^N\eqno{(3.34)}$$
for $z\in Z_1,\,(x',\xi')\in{\rm supp}\,\eta$, and
$$\left|\partial_x^\alpha \left(e^{i\varphi/h}K_2\right)\right|\le
C_{\alpha,N}\left(\frac{h}{\langle\xi'\rangle }\right)^{N-\ell_\alpha}\eqno{(3.35)}$$
otherwise, with possibly a new $\ell_\alpha>0$ independent of $N$. Similar estimates hold for the function $K_1$, too. Indeed, 
observe that on supp$\,\left[\Delta,\Phi_\delta\right]$ we have
$\delta\min\{1,|\rho|^3\}\le x_1\le 2\delta\min\{1,|\rho|^3\}$, and hence
$$\left|e^{i\varphi/h}\right|\le e^{-{\rm Im}\,\varphi/h}\le e^{-x_1{\rm Im}\,\rho/2h}$$
 $$
\le\left\{
\begin{array}{ll}
e^{-C|\rho|^2|{\rm Im}\,z|/h}, &z\in Z_1,\,(x',\xi')\in{\rm supp}\,\eta,\\
e^{-C\langle\xi'\rangle/h}, &\mbox{otherwise},
\end{array}
\right.
\eqno{(3.36)}$$
with some constant $C>0$. Using (3.36) one can easily get that the estimates (3.34) and (3.35) are satisfied with $K_2$ replaced by $K_1$.
Therefore, the function $K$ satisfies the bounds
$$\left|\partial_x^\alpha \left(e^{i\varphi/h}K\right)\right|
\le\left\{
\begin{array}{ll}
C_{\alpha,N}h^{2\epsilon N-\ell_\alpha}, &z\in Z_1(1/2-\epsilon),\,(x',\xi')\in{\rm supp}\,\eta,\\
C_{\alpha,N}\left(\frac{h}{\langle\xi'\rangle }\right)^{N-\ell_\alpha}, &\mbox{otherwise}.
\end{array}
\right.
\eqno{(3.37)}$$
Moreover, since $|\rho|\ge Const>0$ on supp$(1-\chi)$, in the case when $\psi=\psi^0(1-\chi)$ we obtain that (3.37) holds with 
$Z_1(1/2-\epsilon)$ replaced by $Z_1(1-\epsilon)$ and $2\epsilon$ replaced by $\epsilon$. Note now that the kernel, $L_\alpha$, of the operator 
${\cal K}_\alpha$ is given by
$$L_\alpha(x,y')=(2\pi h)^{-d+1}\int e^{\frac{i}{h}\langle y',\xi'\rangle}\partial_x^\alpha\left(e^{\frac{i}{h}\varphi(x,\xi',z))}K(x,\xi',h,z)\right)
d\xi'.$$
If $N$ is taken large enough, (3.37) implies the bounds
$$\left|L_\alpha(x,y')\right|
\le\left\{
\begin{array}{ll}
C_{\alpha,N}h^{2\epsilon N-\ell_\alpha}, &z\in Z_1(1/2-\epsilon),\\
C_{\alpha,N}h^{N-\ell_\alpha}, &z\in Z_2\cup Z_3,
\end{array}
\right.
\eqno{(3.38)}$$
with a new $\ell_\alpha>0$ independent of $N$. When $\psi=\psi^0(1-\chi)$, (3.38) holds with 
$Z_1(1/2-\epsilon)$ replaced by $Z_1(1-\epsilon)$ and $2\epsilon$ replaced by $\epsilon$.
Clearly, (3.27) follows from (3.30) and (3.38).
\eproof

In the case when $\psi=\psi^0$, by (3.23) we have
$$\left|\partial_{x'}^\alpha\partial_{\xi'}^\beta a_{k,j}\right|\le C_{k,j,\alpha,\beta}|\rho|^{-3k-4j-2(|\alpha|+|\beta|)}$$
on supp$\,\eta$, and 
$$\left|\partial_{x'}^\alpha\partial_{\xi'}^\beta a_{k,j}\right|\le C_{k,j,\alpha,\beta}\langle\xi'\rangle^{-j-|\beta|}$$
on supp$\,(1-\eta)$. Since $|\rho|\ge C\sqrt{|{\rm Im}\,z|}$ for $z\in Z_1$, $(x',\xi')\in{\rm supp}\,\chi$ and $|\rho|\ge C>0$ for $z\in Z_2\cup Z_3$ 
or $z\in Z_1$, $(x',\xi')\in{\rm supp}\,(1-\chi)$, we get
$$\left|\partial_{x'}^\alpha\partial_{\xi'}^\beta a_{k,j}\right|\le C_{k,j,\alpha,\beta}|{\rm Im}\,z|^{-3k/2-2j-|\alpha|-|\beta|}$$
for $z\in Z_1$, $(x',\xi')\in{\rm supp}\,\chi$, and 
$$\left|\partial_{x'}^\alpha\partial_{\xi'}^\beta a_{k,j}\right|\le C_{k,j,\alpha,\beta}\langle\xi'\rangle^{-j-|\beta|}$$
otherwise. Hence $(1-\chi)a_{k,j}\in{\cal S}_0^{-j}$, $h^{k+j}a_{k,j}\in{\cal S}_{1/2-\epsilon}^{-j}$ uniformly in $z\in Z_1(1/2-\epsilon)$,
$a_{k,j}\in{\cal S}_0^{-j}$ for $z\in Z_2\cup Z_3$. Therefore, we have
$${\rm Op}_h(\eta a_{k,j})-{\rm Op}_h(\eta_1){\rm Op}_h(\eta a_{k,j})={\cal O}(h^\infty):L^2(\Gamma)\to H_h^m(\Gamma)$$
for every integer $m\ge 0$, where $\eta_1\in C_0^\infty(T^*\Gamma)$ is such that $\eta_1=1$ on supp$\,\eta$.
In view of (2.4) this implies
$$\left\|{\rm Op}_h(\eta a_{k,j})\right\|_{L^2(\Gamma)\to H_h^m(\Gamma)}\le\left\|{\rm Op}_h(\eta a_{k,j})\right\|_{L^2(\Gamma)\to L^2(\Gamma)}
\left\|{\rm Op}_h(\eta_1)\right\|_{L^2(\Gamma)\to H_h^m(\Gamma)}+{\cal O}(h^\infty)$$
$$\le\left\{
\begin{array}{ll}
C_{k,j,m}|{\rm Im}\,z|^{-3k/2-2j}, &z\in Z_1(1/2-\epsilon),\\
C_{k,j,m}, &z\in Z_2\cup Z_3,
\end{array}
\right.
\eqno{(3.39)}
$$
for every integer $m\ge 0$.
In view of Proposition 2.3 we also have
$$\left\|{\rm Op}_h((1-\eta)a_{k,j})\right\|_{L^2(\Gamma)\to H_h^j(\Gamma)}\le C_{k,j}.\eqno{(3.40)}$$
By (3.39) and (3.40) we conclude
$$\left\|{\rm Op}_h(a_{k,j})\right\|_{L^2(\Gamma)\to H_h^j(\Gamma)}\le
\left\{
\begin{array}{ll}
C_{k,j}|{\rm Im}\,z|^{-3k/2-2j}, &z\in Z_1(1/2-\epsilon),\\
C_{k,j}, &z\in Z_2\cup Z_3.
\end{array}
\right.
\eqno{(3.41)}
$$
By Lemma 3.7 we also have
$$\left\|{\rm Op}_h((1-\eta)(a_{k,j}-\widetilde a_{k,j}))\right\|_{L^2(\Gamma)\to H_h^{j+1}(\Gamma)}\le
C_{k,j},\quad z\in Z_1\cup Z_2\cup Z_3.
\eqno{(3.42)}
$$
In the case when $\psi=\psi^0(1-\chi)$, the functions $a_{k,j}$ vanish on supp$\,\chi$, and hence $a_{k,j}\in{\cal S}_0^{-j}$ for
$z\in Z_1$. Therefore, in this case the estimate (3.41) holds with $|{\rm Im}\,z|$ replaced by $1$ and $Z_1(1/2-\epsilon)$
replaced by $Z_1$.

We are ready now to prove Theorem 3.1. If $s=0$ we put $b_0^\psi=-ih(1-\eta)\widetilde a_{1,0}$, and if $s\ge 1$ we put
$$b_s^\psi=-i\sum_{j=0}^{s-1}h^{j+1}a_{1,j}-c_sh^s\rho^{-s-1}zn_s\psi-ih^{s+1}(1-\eta)\widetilde a_{1,s}.$$
In view of Lemma 3.6, the function $b_s^\psi$ is independent of all $n_\ell$ with $\ell\ge s$. If we take $N$ big enough, we can decompose the function
$\tau_\psi$ as
$$\tau_\psi=\rho\psi+b_s^\psi+c_sh^s\rho^{-s-1}zn_s\psi+\widetilde b_s^\psi$$
where
$$\widetilde b_s^\psi=-ih^{s+1}\eta a_{1,s}-ih^{s+1}(1-\eta)(a_{1,s}-\widetilde a_{1,s})-i\sum_{j=s+1}^{N-1}h^{j+1}a_{1,j}.$$
By (3.39), (3.41) and (3.42) we have
$$\left\|{\rm Op}_h(\widetilde b_s^\psi)\right\|_{L^2(\Gamma)\to H_h^{s+1}(\Gamma)}\le
\left\{
\begin{array}{ll}
C_sh^{s+1}|{\rm Im}\,z|^{-3/2-2s}, &z\in Z_1(1/2-\epsilon),\\
C_sh^{s+1}, &z\in Z_2\cup Z_3.
\end{array}
\right.
\eqno{(3.43)}
$$
Moreover, if $\psi=\psi^0(1-\chi)$, the estimate (3.43) holds with $|{\rm Im}\,z|$ replaced by $1$ and 
$Z_1(1/2-\epsilon)$ replaced by $Z_1$. We would like to apply Lemma 3.8 with $m=s+1$. To this end we take $N$ big enough to arrange that
$$\epsilon N-\ell_{s+1}>s+1.$$
By (3.27) and (3.43) we get
$$\left\|{\cal N}(h,z){\rm Op}_h(\psi)-{\rm Op}_h(\rho\psi+b_s^\psi+c_sh^s\rho^{-s-1}zn_s\psi)\right\|_{L^2(\Gamma)\to H_h^{s+1}(\Gamma)}$$
$$\le
\left\{
\begin{array}{ll}
C_sh^{s+1}|{\rm Im}\,z|^{-3/2-2s}, &z\in Z_1(1/2-\epsilon),\\
C_sh^{s+1}, &z\in Z_2\cup Z_3,
\end{array}
\right.
\eqno{(3.44)}
$$
if $\psi=\psi^0$. 
Moreover, if $\psi=\psi^0(1-\chi)$, the estimate (3.44) holds with $|{\rm Im}\,z|$ replaced by $1$ and 
$Z_1(1/2-\epsilon)$ replaced by $Z_1(1-\epsilon)$. 

We will now use a partition of the unity on $\Gamma$. We can find functions $\{\psi_j^0\}_{j=1}^J$ such that $\sum_{j=1}^J\psi_j^0=1$
and (3.44) is valid with $\psi$ replaced by each $\psi_j$, where $\psi_j$ is defined by replacing in the definiton of $\psi$ the function $\psi^0$ by $\psi_j^0$. 
Summing up all the estimates we get (3.8) and (3.9), respectively. 
\eproof

\section{Improved estimates}

To prove Theorem 1.1 we actually need the following improved version of Theorem 3.1.

\begin{Theorem} Let $z\in Z_1(1/2-\epsilon)$, $0<\epsilon\ll 1$. Then, for every integer $s\ge 1$ there are an operator ${\cal B}_s$ independent of all $n_k$ with $k\ge s$ and an operator 
$${\cal A}_s={\cal O}_s(h^{-s}):H_h^{s+1}(\Gamma)\to 
L^2(\Gamma)\eqno{(4.1)}$$
independent of all $n_k$ with $k\ge 1$ such that
$$\left\|{\cal A}_s{\cal N}(h,z)-{\cal B}_s-n_sI\right\|_{L^2(\Gamma)\to L^2(\Gamma)}\le 
C_sh|{\rm Im}\,z|^{-3s/2-1}\eqno{(4.2)}$$
where $I$ denotes the identity. 
If $z\in Z_2\cup Z_3$, then (4.2) holds with $|{\rm Im}\,z|$ replaced by $1$.
\end{Theorem}

{\it Proof.} Recall that by (3.5), (3.6), (3.7), we have
that for every integer $k$, $h^{\frac{k_-}{4}}\rho^k\in {\cal S}_{1/2-\epsilon}^k$ uniformly in $z\in Z_1(1/2-\epsilon)$ and 
$\rho^k\in {\cal S}_0^k$ if $z\in Z_2\cup Z_3$. We would like to apply Proposition 2.2 with
$$a^+_s=\left(c_sh^s\rho^{-s-1}z\right)^{-1},\quad a^-_s=c_sh^s\rho^{-s-1}z.$$
Using (3.2) one can easily check that (2.5) is satisfied with $\mu=h|{\rm Im}\,z|^{-2}={\cal O}(h^{2\epsilon})$. By (2.6) we get
$$\left\|{\rm Op}_h(a^+_s){\rm Op}_h(a^-_s)-I\right\|_{L^2(\Gamma)\to L^2(\Gamma)}\le Ch^{2\epsilon}\le 1/2\eqno{(4.3)}$$
if $h$ is taken small enough. It follows from (4.3) that the operator ${\rm Op}_h(a^-_s)$ is invertible with an inverse
$${\cal A}_s:=\left({\rm Op}_h(a^-_s)\right)^{-1}=\left({\rm Op}_h(a^+_s){\rm Op}_h(a^-_s)\right)^{-1}{\rm Op}_h(a^+_s).$$
Since $h^sa^+_s\in {\cal S}_{1/2-\epsilon}^{s+1}$ uniformly in $z$, by Proposition 2.3 we have  
$$h^s{\rm Op}_h(a^+_s)={\cal O}_s(1):H_h^{s+1}(\Gamma)\to 
L^2(\Gamma)$$ 
which implies (4.1). 
By (3.27) and (4.1),
$$\left\|{\cal A}_s{\cal N}(h,z)-{\cal A}_s{\cal T}(h,z)\right\|_{L^2(\Gamma)\to L^2(\Gamma)}\le 
C_{N,s}h^{\epsilon N-s-\ell_{s+1}}\le h\eqno{(4.4)}$$
if $N$ is taken large enough, where ${\cal T}=\sum_{j=1}^J{\cal T}_{\psi_j}$. On the other hand, we can write
$${\cal A}_s{\cal T}={\cal B}_s+n_sI+\widetilde{\cal B}_s$$
where
$${\cal B}_s={\cal A}_s{\rm Op}_h(\rho+b_s),\quad b_s=\sum_{j=1}^Jb_s^{\psi_j},$$
$$\widetilde{\cal B}_s={\cal A}_s{\rm Op}_h(\widetilde b_s),\quad \widetilde b_s=\sum_{j=1}^J\widetilde b_s^{\psi_j}.$$
Clearly, the operator ${\cal B}_s$ is independent of all $n_k$ with $k\ge s$ because so is the function $b_s$. Therefore, it follows from (4.4) that 
to prove (4.2) it suffices to prove the bound
 $$\left\|{\rm Op}_h(a^+_s){\rm Op}_h(\widetilde b_s)\right\|_{L^2(\Gamma)\to L^2(\Gamma)}\le 
 \left\{
\begin{array}{ll}
C_sh|{\rm Im}\,z|^{-3s/2-1}, &z\in Z_1(1/2-\epsilon),\\
C_sh, &z\in Z_2\cup Z_3.
\end{array}
\right.
\eqno{(4.5)}$$
In view of Lemmas 3.5 and 3.7, we have $\widetilde b_s=h^{s+1}g_s$ with $g_s\in S_{2,2}^{-3-4s}(|\rho|)+S_{0,1}^{-s-1}(|\rho|)$ uniformly in $h$
as long as $|\rho|^4\ge h$.
Thus, (4.5) is equivalent to
$$\left\|{\rm Op}_h(\rho^{s+1}){\rm Op}_h(g_s)\right\|_{L^2(\Gamma)\to L^2(\Gamma)}\le 
 \left\{
\begin{array}{ll}
C_s|{\rm Im}\,z|^{-3s/2-1}, &z\in Z_1(1/2-\epsilon),\\
C_s, &z\in Z_2\cup Z_3.
\end{array}
\right.
\eqno{(4.6)}$$
To prove (4.6) observe that $\rho^{s+1}g_s\in S_{2,2}^{-2-3s}(|\rho|)+S_{0,1}^0(|\rho|)$ uniformly in $h$, which yields the bounds
$$\left|\partial_{x'}^\alpha\partial_{\xi'}^\beta\left(\rho^{s+1}g_s\right)\right|\le 
 \left\{
\begin{array}{ll}
C_{s,\alpha,\beta}|{\rm Im}\,z|^{-3s/2-1-|\alpha|-|\beta|}, &z\in Z_1(1/2-\epsilon),\\
C_{s,\alpha,\beta}, &z\in Z_2\cup Z_3.
\end{array}
\right.
\eqno{(4.7)}$$
By (2.4) and (4.7) we get
$$\left\|{\rm Op}_h(\rho^{s+1}g_s)\right\|_{L^2(\Gamma)\to L^2(\Gamma)}\le 
 \left\{
\begin{array}{ll}
C_s|{\rm Im}\,z|^{-3s/2-1}, &z\in Z_1(1/2-\epsilon),\\
C_s, &z\in Z_2\cup Z_3.
\end{array}
\right.
\eqno{(4.8)}$$
On the other hand, applying Proposition 2.2 with $a^+=\rho^{s+1}$ and $a^-=g_s$ yields the bound
$$\left\|{\rm Op}_h(\rho^{s+1}){\rm Op}_h(g_s)-{\rm Op}_h(\rho^{s+1}g_s)\right\|_{L^2(\Gamma)\to L^2(\Gamma)}$$
$$\le 
 \left\{
\begin{array}{ll}
C_sh|{\rm Im}\,z|^{-3s/2-3}, &z\in Z_1(1/2-\epsilon),\\
C_sh, &z\in Z_2\cup Z_3.
\end{array}
\right.
\eqno{(4.9)}$$
Clearly, (4.6) follows from (4.8) and (4.9).
\eproof

\section{Proof of Theorem 1.1}

Define the DN maps ${\cal N}_j(\lambda)$, $j=1,2$, by
$${\cal N}_j(\lambda)f=\partial_\nu u_j|_\Gamma$$
where $\nu$ is the Euclidean unit inner normal to $\Gamma$ and $u_j$ is the solution to the equation
$$\left\{
\begin{array}{lll}
 \left(\Delta+\lambda^2n_j(x)\right)u_j=0&\mbox{in}& \Omega,\\
 u_j=f&\mbox{on}&\Gamma,
\end{array}
\right.
\eqno{(5.1)}
$$
and consider the operator
$$T(\lambda)={\cal N}_1(\lambda)-{\cal N}_2(\lambda).$$
Clearly, $\lambda$ is a transmission eigenvalue if there exists a non-trivial function $f$ such that $T(\lambda)f=0$.
Thus Theorem 1.1 is a consequence of the following

\begin{Theorem} Under the conditions of Theorem 1.1, 
the operator $T(\lambda)$ sends $L^2(\Gamma)$ into $H^{j+1}(\Gamma)$. 
Moreover, there exists a constant $C>0$ such that $T(\lambda)$ is invertible for 
$|{\rm Im}\,\lambda|\ge C({\rm Re}\,\lambda+1)^{1-\kappa_j}$ with an inverse satisfying in this region the bound
$$\left\|T(\lambda)^{-1}\right\|_{H^{j+1}(\Gamma)\to L^2(\Gamma)}\lesssim|\lambda|^{j-1}\eqno{(5.2)}$$
where the Sobolev space is equipped with the classical norm.
\end{Theorem}

{\it Proof.} We make our problem semi-classical by putting $h=|{\rm Re}\,\lambda^2|^{-1/2}$, 
$z=h^2\lambda^2=\pm 1+i{\rm Im}\,z$,
if $|{\rm Re}\,\lambda^2|\ge|{\rm Im}\,\lambda^2|$, $\pm {\rm Re}\,\lambda^2>0$, 
and $h=|{\rm Im}\,\lambda^2|^{-1/2}$, $z=h^2\lambda^2={\rm Re}\,z+i$,
if $|{\rm Re}\,\lambda^2|\le|{\rm Im}\,\lambda^2|$. Clearly, $h\sim|\lambda|^{-1}$. We set 
${\cal N}_j(h,z)=-ih{\cal N}_j(\lambda)$ and
$$T(h,z)={\cal N}_1(h,z)-{\cal N}_2(h,z).$$
We now apply Theorem 4.1 with $s=j\ge 1$. In view of the conditions (1.5) and (1.6), we get
$$\left\|{\cal A}_jT(h,z)-(n^{(1)}_j-n^{(2)}_j)I\right\|_{L^2(\Gamma)\to L^2(\Gamma)}\le 
C_jh|{\rm Im}\,z|^{-3j/2-1}\eqno{(5.3)}$$
for $z\in Z_1(1/2-\epsilon)$, where
$$n^{(\ell)}_j=(j!)^{-1}\partial_\nu^j n_\ell|_\Gamma,\quad \ell=1,2.$$
 When $z\in Z_2\cup Z_3$, the estimate (5.3) holds
with $|{\rm Im}\,z|$ replaced by $1$.
 It follows from (5.3) that the operator $(n^{(1)}_j-n^{(2)}_j)^{-1}{\cal A}_jT(h,z)$ is invertible for
$z\in Z_1(1/2-\epsilon)$, $|{\rm Im}\,z|\ge (C'_jh)^{1/(3j/2+1)}$, and for $z\in Z_2\cup Z_3$, $h$ small enough.
Hence so is $T(h,z)$ and we have the bound
$$\left\|T(h,z)^{-1}\right\|_{H_h^{j+1}(\Gamma)\to L^2(\Gamma)}\le {\cal O}(1)\left\|{\cal A}_j\right\|_{H_h^{j+1}(\Gamma)\to L^2(\Gamma)}\le
{\cal O}(h^{-j}).\eqno{(5.4)}$$
Now (5.2) follows from (5.4) after passing from $(h,z)$ to $\lambda$ and using the fact that the semi-classical norm 
in $H_h^{j+1}(\Gamma)$ is bounded from above by the classical norm 
in $H^{j+1}(\Gamma)$.
\eproof

It is worth noticing that it follows from the estimate (3.9) that the operator $T(h,z)$ can be inverted outside the glancing region
for much smaller $|{\rm Im}\,z|$. In other words, to improve the eigenvalue-free region (1.7) one has to improve the parametrix in
the glancing region, only. More precisely, we have the following

\begin{prop} Let $z\in Z_1(1-\epsilon)$. Then, under the conditions of Theorem 1.1, there exists an operator
$$\widetilde{\cal A}_j={\cal O}(h^{-j}):H_h^{j+1}(\Gamma)\to 
L^2(\Gamma)$$
such that
$$\left\|T(h,z)\widetilde{\cal A}_j-{\rm Op}_h(1-\chi)\right\|_{H_h^{j+1}(\Gamma)\to H_h^{j+1}(\Gamma)}\le 
Ch.\eqno{(5.5)}$$
When $z\in Z_2\cup Z_3$, the estimate (5.5) holds with $\chi$ replaced by $0$.
\end{prop}

{\it Proof.} By (3.9) with $s=j$ we have 
$$\left\|T(h,z){\rm Op}_h(1-\chi)-{\rm Op}_h\left((1-\chi)c_jh^j\rho^{-j-1}z(n^{(1)}_j-n^{(2)}_j)\right)\right\|_{L^2(\Gamma)\to H_h^{j+1}(\Gamma)}\le 
C_jh^{j+1}\eqno{(5.6)}$$
for $z\in Z_1(1-\epsilon)$. Let $\chi_1\in C_0^\infty(T^*\Gamma)$ be such that $\chi=1$ 
on supp$\,\chi_1$, $\chi_1=1$ in a neighbourhood of $\Sigma$, and set 
$$\widetilde a_j^+=(1-\chi)\rho^{-j-1}c_jz(n^{(1)}_j-n^{(2)}_j),\quad \widetilde a_j^-=(1-\chi_1)\rho^{j+1}\left(c_jz(n^{(1)}_j-n^{(2)}_j)\right)^{-1}.$$
We have $\widetilde a_j^+\in {\cal S}_0^{-j-1}$, $\widetilde a_j^-\in {\cal S}_0^{j+1}$ and $\widetilde a_j^+\widetilde a_j^-=1-\chi$. We now apply Proposition 2.4 with $\widetilde a_j^+$
and $\widetilde a_j^-$ in place of $a^+$ and $a^-$.  We have
$$\left\|{\rm Op}_h(\widetilde a_j^+){\rm Op}_h(\widetilde a_j^-)-{\rm Op}_h(1-\chi)\right\|_{H_h^{j+1}(\Gamma)\to H_h^{j+1}(\Gamma)}\le Ch.\eqno{(5.7)}$$
Clearly, (5.5) follows from (5.6) and (5.7) with $\widetilde{\cal A}_j=h^{-j}{\rm Op}_h(1-\chi){\rm Op}_h(\widetilde a_j^-)$.
\eproof

\end{document}